\documentclass[11pt,reqno]{amsart}
\usepackage{amsfonts,amsmath,amssymb}
\newtheorem{Theorem}{Theorem}

\newtheorem{Lemma}{Lemma}

\def\beq#1#2\eeq{%
        \begin{equation}%
        \label{#1}%
            #2%
        \end{equation}%
    }

\usepackage{color}
\usepackage{graphicx}
\usepackage{epstopdf}

\theoremstyle{plain}

\theoremstyle{remark}

\theoremstyle{definition}

\title[Coway and Arnold]{Conway river and Arnold sail}

\author{ K. Spalding}\address{Department of Mathematical Sciences,
Loughborough University, Loughborough LE11 3TU, UK}
\email{K.Spalding@lboro.ac.uk}

\author{A.P. Veselov}
\address{Department of Mathematical Sciences,
Loughborough University, Loughborough LE11 3TU, UK  and Moscow State University, Moscow 119899, Russia}
\email{A.P.Veselov@lboro.ac.uk}

\begin{document}

\maketitle

\centerline{\it To the memory of Vladimir I. Arnold, who would be 80 now}

\begin{abstract}
We establish a simple relation between two geometric constructions in number theory: the Conway river of a real indefinite binary quadratic form and the Arnold sail of the corresponding pair of lines.
\end{abstract}

\section{Introduction}

In 1895 Felix Klein proposed the following geometric representation of continued fractions. For an irrational real number $\omega$ consider the ray $y = \omega x$ on the plane with the integer lattice. Let us quote Klein (1924):

\begin{quote}
Imagine pegs or needles affixed at all the integral points, and wrap a tightly drawn string about the sets of pegs to the right and to the left of the $\omega$-ray, then the vertices of the two convex strong-polygons which bound our two point sets will be precisely the points $(p_{\nu}, q_{\nu})$ whose coordinates are the numerators and denominators of the successive convergents to $\omega$, the left polygon having the even convergents, the right one the odd. This gives a new and, one may well say, an extremely graphic definition of a continued fraction.
\end{quote}

\begin{figure}[h]
	\centering
	\includegraphics[scale=0.5]{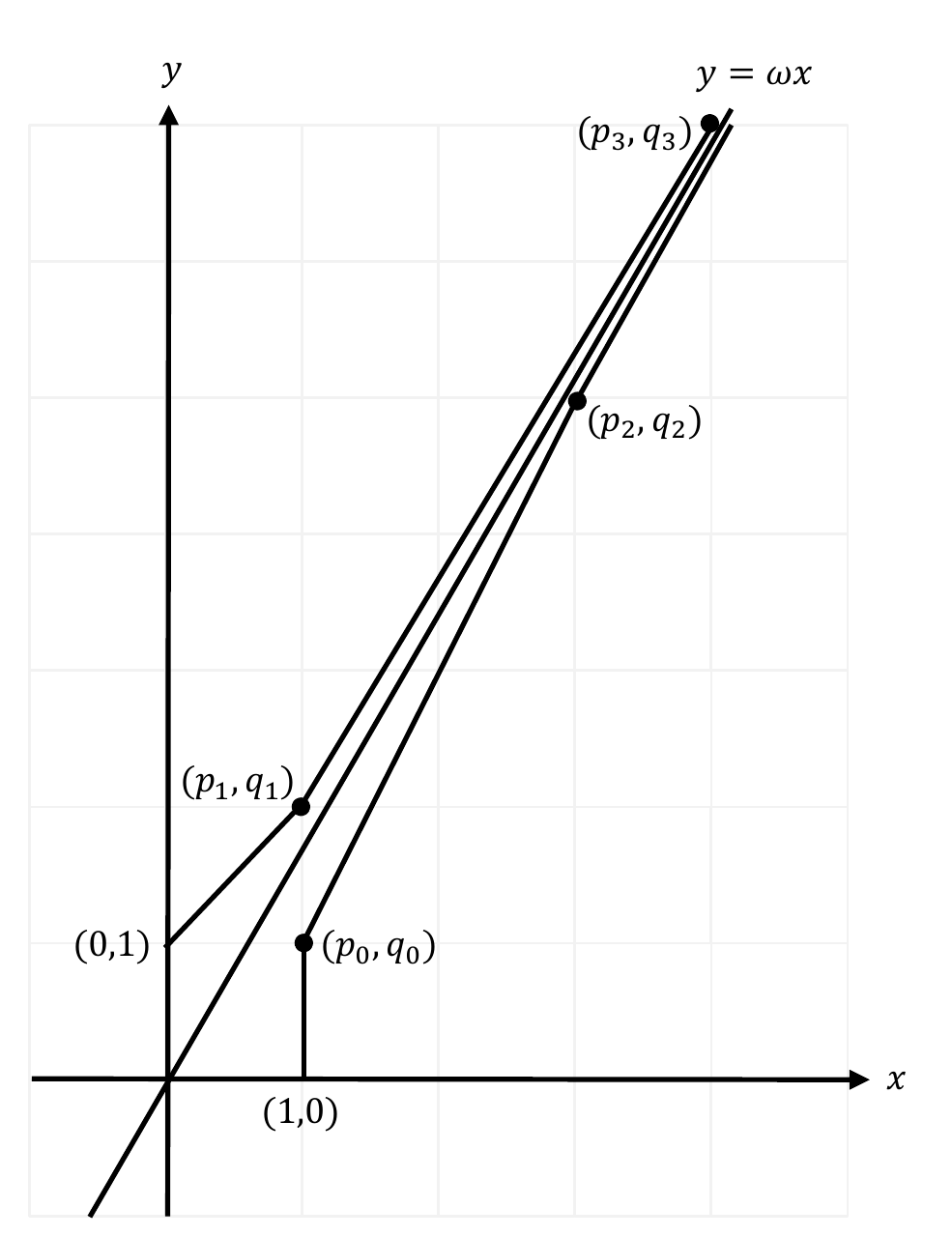} \hspace{8pt} \includegraphics[scale=0.5]{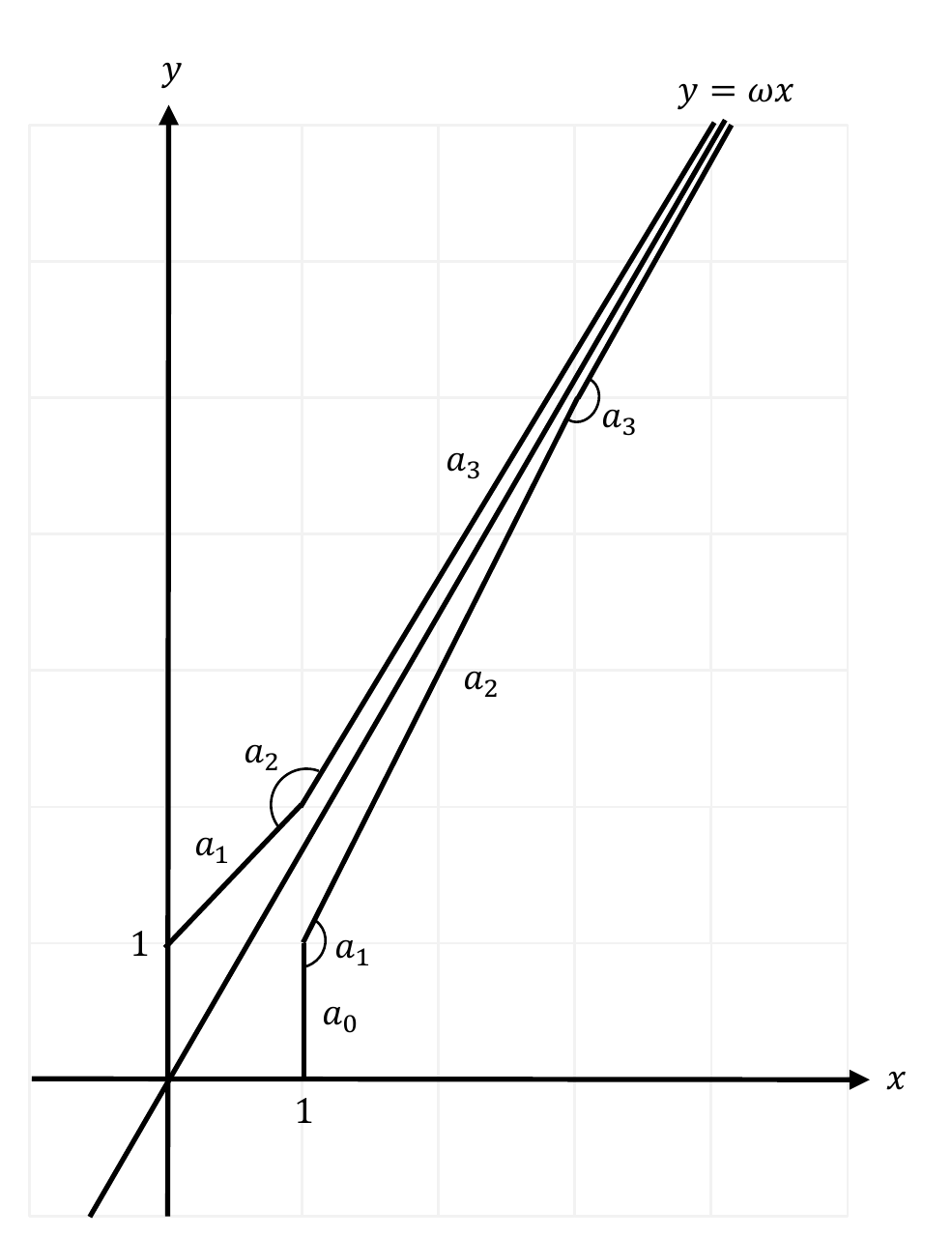}
	\caption{Klein's construction and the corresponding LLS sequence} \label{fig:Klein}
	\end{figure}

Many years later Vladimir I. Arnold (1998) revitalised this point of view, mainly with an emphasis on multi-dimensional generalisations. In particular, for a polyhedral cone he introduced the notion of the {\it sail} as the boundary of the convex hull of the integer points inside it.  In dimension 2 the sail of the angle formed by the $\omega$-ray and $x$-axis is precisely Klein's construction of the continued fraction expansion of $\omega$ (see Figure \ref{fig:Klein}). 

\begin{figure}[h]
	\centering
	\includegraphics[scale=0.2]{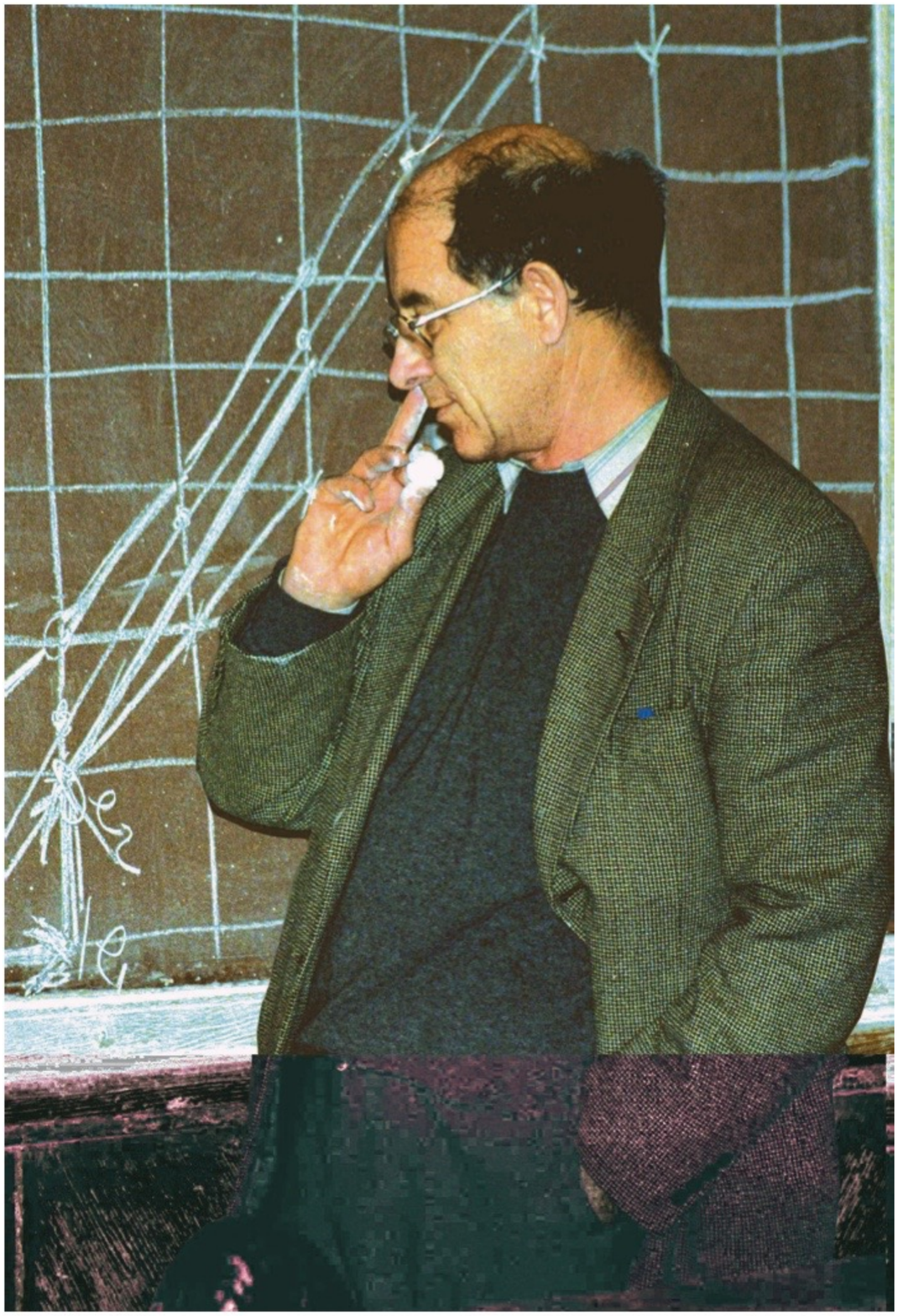} \hspace{8pt}  \includegraphics[height=56mm]{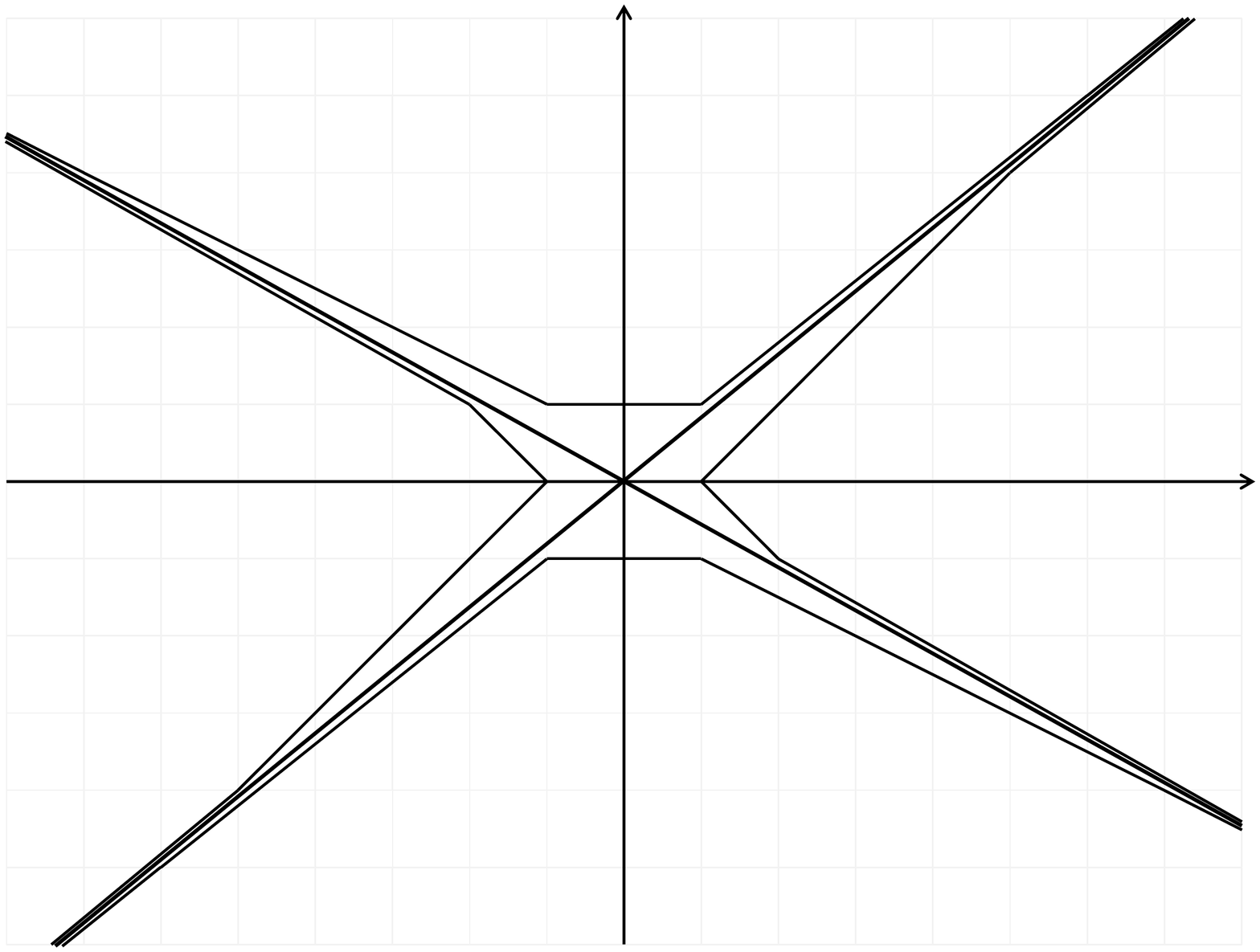}
	\caption{Arnold and the sails for a pair of lines.} \label{fig:arnold_4graph}
	\end{figure}

This line was developed in more detail by Karpenkov (2013), who, importantly for us, introduced the \emph{lattice length sine (LLS) sequence} of positive integers $(a_i), \, i \in \mathbb Z$ of the sail and proved a remarkable edge-angle duality between the sails of the adjacent angles (see Figure \ref{fig:Klein}).  He also linked this with the theory of indefinite binary quadratic forms. Indeed, the zero set of such a form is a pair of lines, forming 4 angles with 4 sails, which are either isomorphic or dual to each other (see Figure \ref{fig:arnold_4graph}).\footnote {Photo of V.I. Arnold is reproduced from (Arnold 2002) with permission from MCCME.}

\begin{figure}[h]
\begin{center}
 \includegraphics[height=40mm]{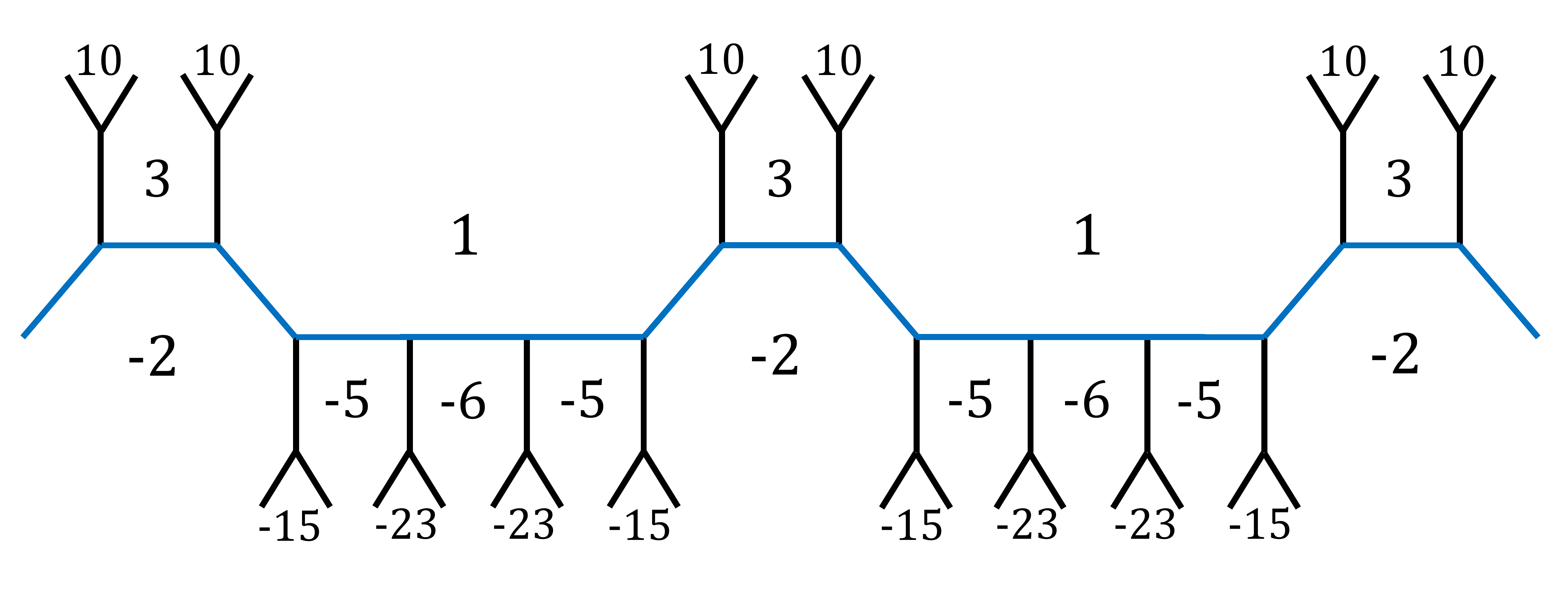}  
\caption{\small Conway river for the quadratic form  $Q=x^2-2xy-5y^2.$} \label{fig:river}
\end{center}
\end{figure}

On the other hand, John H. Conway (1997) proposed the notion of the {\it topograph} of a binary quadratic form $Q$ as a graphical way to visualise the values of $Q$ on a planar binary tree (see section 3 below). For indefinite quadratic forms he introduced the notion of the {\it river}, which is a path on the topograph separating positive and negative values of $Q$ (see Figure \ref{fig:river}).

The main result of this paper is the following simple relation between the Conway river and the corresponding Arnold sail.  

\begin{Theorem}
Let $Q(x,y)$ be a real indefinite binary quadratic form and consider the Arnold sail of the pair of lines given by $Q(x,y)=0,$ assuming that the origin is the only integer point on them.

Then the corresponding LLS sequence $(a_i), \, i \in \mathbb Z$ coincides with the sequence of the left- and right-turns of the Conway river on the topograph of $Q.$ This determines the river uniquely up to the action of the group $PGL(2, \mathbb Z)$ on the topograph and a change of direction.
\end{Theorem}

For example, for $Q=x^2-2xy-5y^2$ one can check that the corresponding LLS sequence is periodic, equal to $\dots 4,2,4,2,4,2, \dots$, which is exactly the sequence of left-right turns $\dots LLLLRRLLLLRR\dots$ of the (properly oriented) Conway river in Figure \ref{fig:river}. 

The proof is simple and essentially follows from the results of Karpenkov (2013) combined with more detailed analysis of the Conway river from (Spalding and Veselov 2017).

\section{Arnold sail and the LLS sequence of the angles}

We follow here Karpenkov (2013) (see, in particular, Chapters 2 and 4).

Let $A,B,C$ be three distinct integer lattice points on the plane and $\angle ABC$ be the corresponding angle. Define the {\it integer length} $\mathrm{l}l(AB)$ of the segment $AB$ as the number of integer points in the interior of $AB$ plus one, and the {\it integer area} $\mathrm{l}S(\triangle ABC)$ of the triangle $\triangle ABC$ as the index of the sublattice generated by the integer vectors $AB$ and $BC$ in the integer lattice.

The {\it integer sine} of the angle $\angle ABC$ is defined as
\begin{equation}
\label{lsin}
\mathrm{l}\sin \angle ABC=\frac{\mathrm{l}S(\triangle ABC)}{\mathrm{l}l(AB)\mathrm{l}l(BC)}=\frac{|\det(AB,BC)|}{\mathrm{l}l(AB)\mathrm{l}l(BC)}.
\end{equation}
One can check that it is a positive integer and depends only on the angle, and not on the choice of $B$ and $C.$

Consider now a pair of lines given by $y=\alpha x$ and $y=\beta x$ and one of the angles $\angle {\mathcal A} O {\mathcal  B}$ formed by them. 
Let us assume that $\alpha$ and $\beta$ are irrational, so that the origin $O$ is the only integer point on them, and consider the convex hull of the integer points inside $\angle {\mathcal A} O {\mathcal B}$ (excluding $O$). Its boundary is an infinite broken line called the {\it Arnold sail} of the angle $\angle {\mathcal A} O {\mathcal  B}$. 

Let $(A_i), \, i\in \mathbb Z$ be the sequence of vertices of this sail. 
Karpenkov (2013) introduced the following key notion of the \emph{LLS (lattice length sine) sequence} $(a_i), \, i \in \mathbb Z$ of the angle $\angle {\mathcal A} O {\mathcal  B}$ as
\begin{equation}
\label{LLS}
a_{2k} = \mathrm{l}l A_k A_{k+1}, \quad
a_{2k-1} = \mathrm{l} \sin \left(\angle A_{k-1} A_k A_{k+1} \right).
\end{equation}
Karpenkov proved that the LLS sequence determines the angle uniquely up to an integer affine transformation. Note that the sequence is defined up to a shift of indices and depends on the orientation of the sail. 

When the angle is formed by the $x$-axis and $\omega$-ray the corresponding LLS sequence is semi-infinite and gives precisely the continued fraction representation of $\omega$ (see Figure \ref{fig:Klein}): 
$$\omega=[a_0, a_1, a_2, \dots]:=a_0+\frac{1}{a_1 +\frac{1}{a_2+ \dots}}.$$

Let us look at the sail of the angle formed by the $\omega$-ray and $y$-axis.
Let $B_0B_1B_2 \ldots$ be the sequence of vertices of the corresponding sail. Then we have the remarkable {\it edge-angle duality} (Karpenkov 2013):
\begin{equation}
\mathrm{l}\sin \left( \angle A_i A_{i+1} A_{i+2} \right) = \mathrm{l} l \left(B_i B_{i+1} \right),
\end{equation}
\begin{equation}
\mathrm{l}\sin \left( \angle B_i B_{i+1} B_{i+2} \right) = \mathrm{l} l \left(A_{i+1} A_{i+2} \right).
\end{equation}
This explains why we do not need to consider the second sail to extract the full continued fraction expansion.
Note that the coordinates of $A_i=(p_{2i}, q_{2i})$ and $B_i=(p_{2i-1}, q_{2i-1})$ are the corresponding denominators and numerators of the continued fraction convergents for $\omega$ (see (Klein 1924) and Figure \ref{fig:Klein}).

For general lines given by $y=\alpha x$ and $y=\beta x$ the corresponding (infinite in both directions) LLS sequence can be considered as {\it a joint continued fraction expansion} of the numbers $\alpha$ and $\beta$ and is related to the rational approximation of the arrangement of these two lines (or, equivalently to the corresponding quadratic form $Q=(y-\alpha x)(y-\beta x)$, see Chapter 10 in (Karpenkov 2013)).

\section{Topograph of binary quadratic form and Conway river}

We follow here the original approach of Conway (1997).

Conway proposed the following nice way to visualise the values of a binary quadratic form
\begin{equation}
\label{Q}
Q(x, y) = ax^2 + hxy + by^2, \quad (x,y) \in \mathbb{Z}^2.
\end{equation}
He considered the case when all the coefficients $a,b,h$ are integer, but his construction works for real coefficients as well.

\begin{figure}[h]
\begin{center}
 \includegraphics[height=41mm]{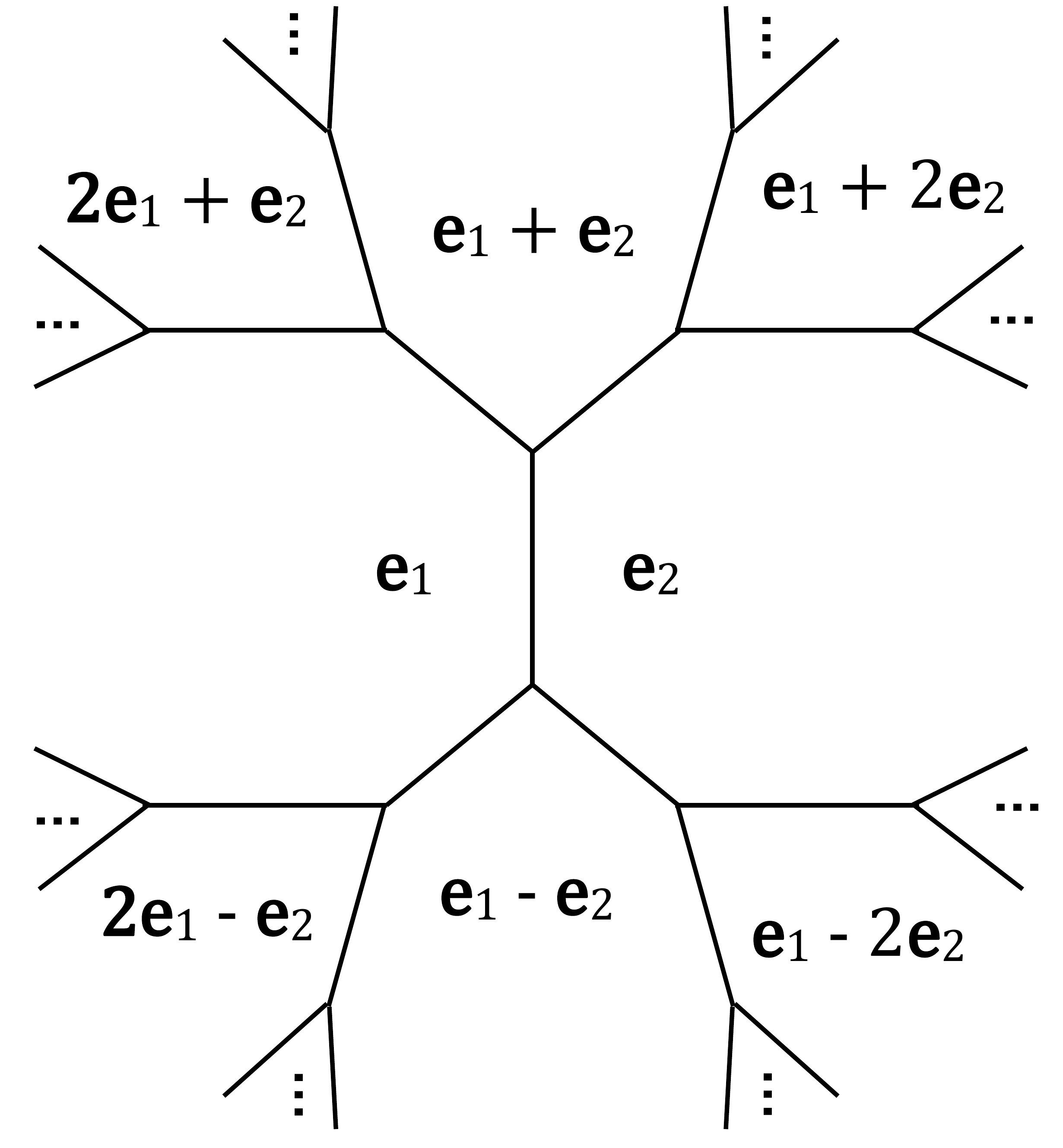}  \hspace{15pt} \includegraphics[scale=0.23]{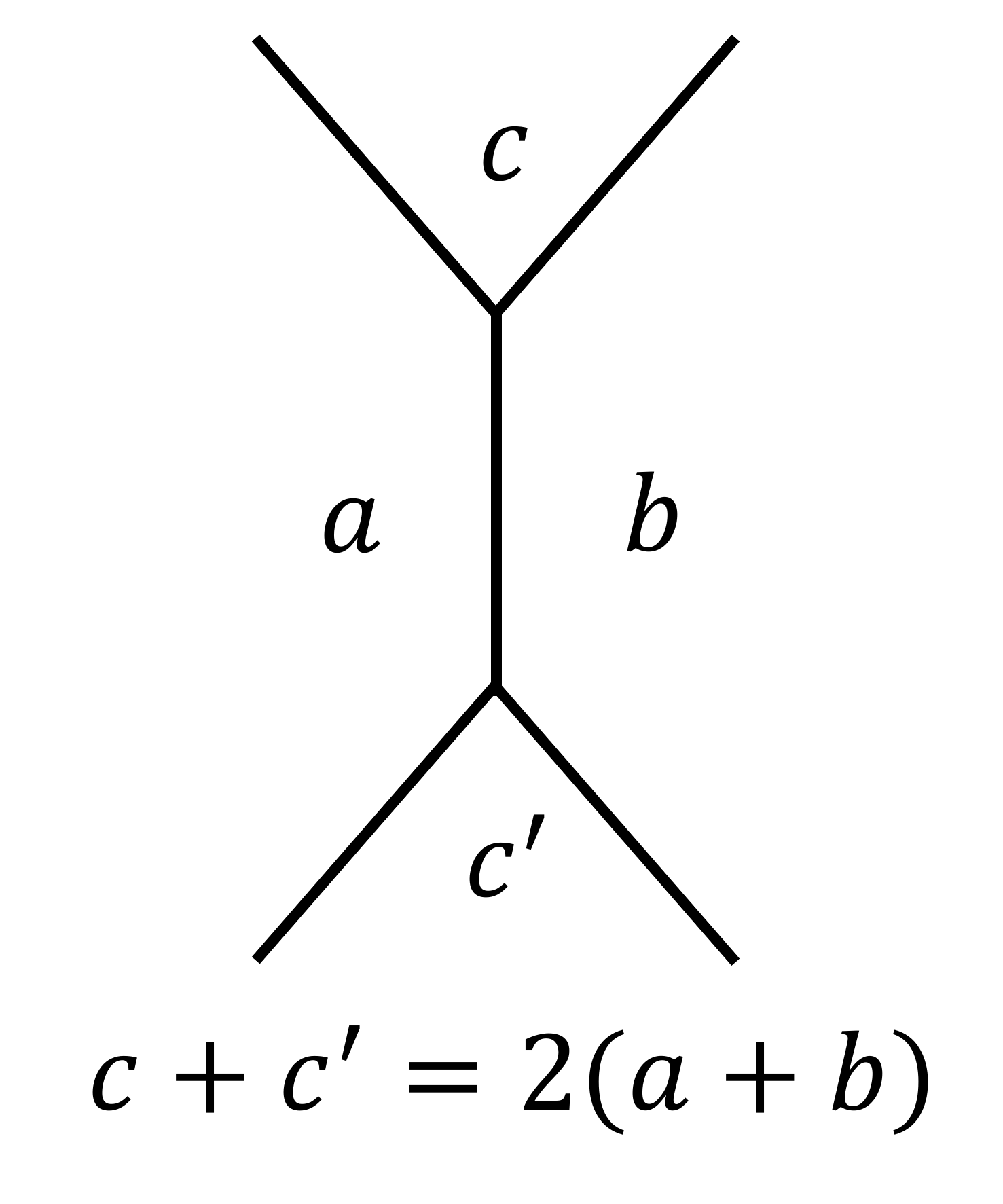}
\caption{\small The superbase tree and arithmetic progression rule for values of quadratic forms.} \label{fig:tree}
\end{center}
\end{figure}

Conway introduced the notions of the {\it lax} vector as a pair $(\pm v), v \in \mathbb Z^2$, and of the \emph{superbase} of the integer lattice $\mathbb Z^2$ as a triple of lax vectors $(\pm e_1, \pm e_2, \pm e_3)$ such that $(e_1, e_2)$ is a basis of the lattice and
		\begin{equation*}
		e_1 + e_2 + e_3 =  0.
		\end{equation*}
It is easy to see that every basis can be included in exactly two superbases, which we can represent using the binary tree embedded in the plane (see Figure \ref{fig:tree}).
The lax vectors live in the complement to the tree (we show only one representative of them), while the superbases correspond to the vertices.
Note that all {\it primitive} lattice vectors, i.e. those which are not multiples of any other lattice vectors, appear on this tree.
		
By taking the values of the form $Q$ on the vectors of the superbase tree, we get what Conway called the {\it topograph} of $Q.$ The idea is to get the values of $Q$ on all primitive lattice vectors in this way.

In particular, if $e_1=(1,0), e_2=(0,1), e_3=-(1,1)$ we have the values 
$Q(e_1)=a, \,\, Q(e_2)=b, \,\, Q(e_3)=c:=a+b+h.$
One can construct the topograph of $Q$ starting from this triple using the {\it arithmetic progression rule} (known in geometry as the {\it parallelogram rule}):
\begin{equation}
\label{apr}
Q(u+v)+Q(u-v)=2(Q(u)+Q(v)), \quad u,v \in \mathbb R^2.
\end{equation}
%
	
We also need to construct the {\it Farey tree} by replacing $v=(p,q)$ on the superbase tree by the corresponding fraction $\frac{p}{q}$ (so that addition of vectors corresponds to taking the Farey mediant of fractions). 

Using the Farey tree, we can label any semi-infinite path $\gamma$ on the tree by a real number $\xi$ such that the limit of the Farey fractions along $\gamma_\xi$ is $\xi$. 

\begin{figure}[h]
\begin{center}
 \includegraphics[height=45mm]{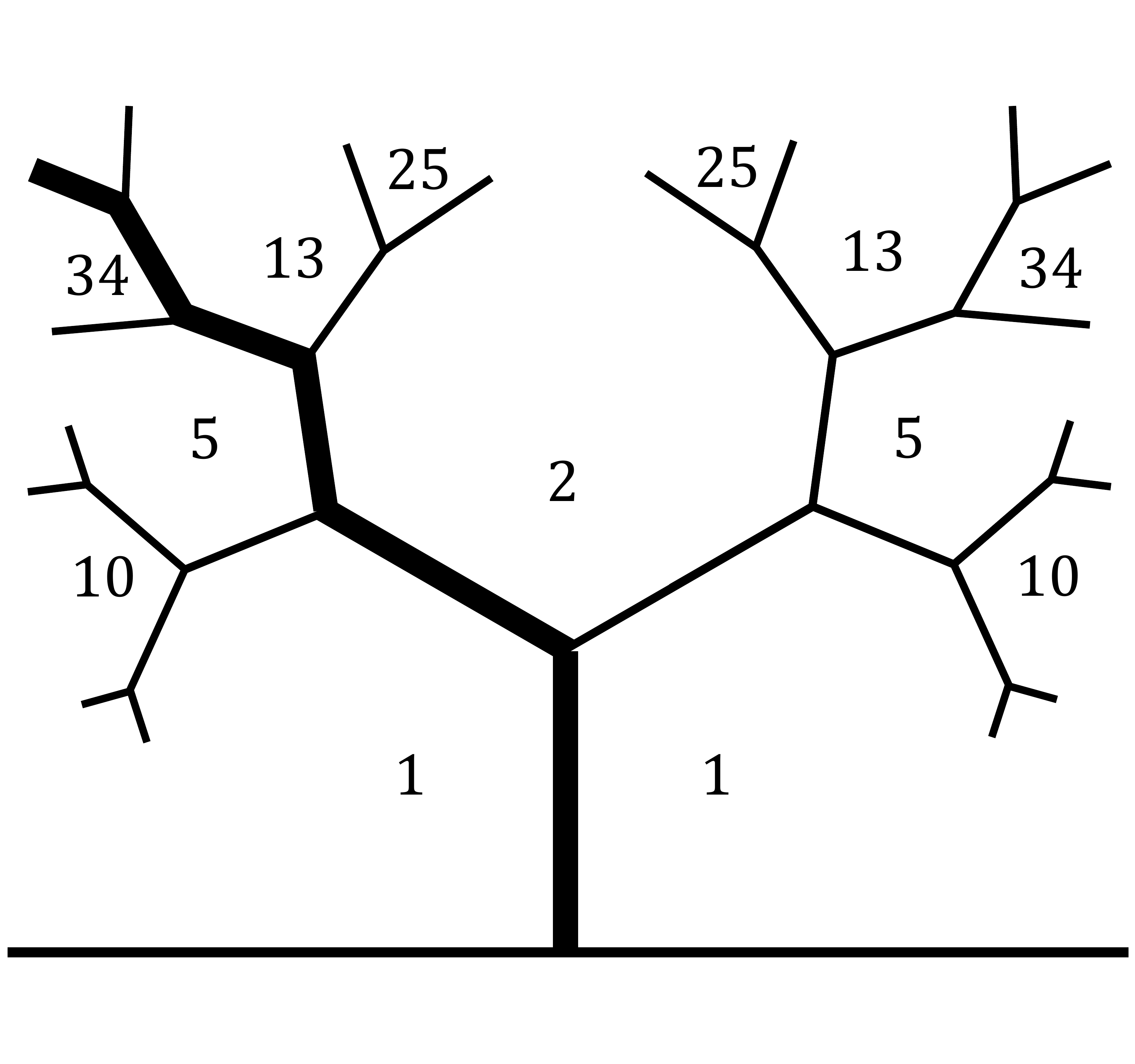}  \hspace{8pt}  \includegraphics[height=45mm]{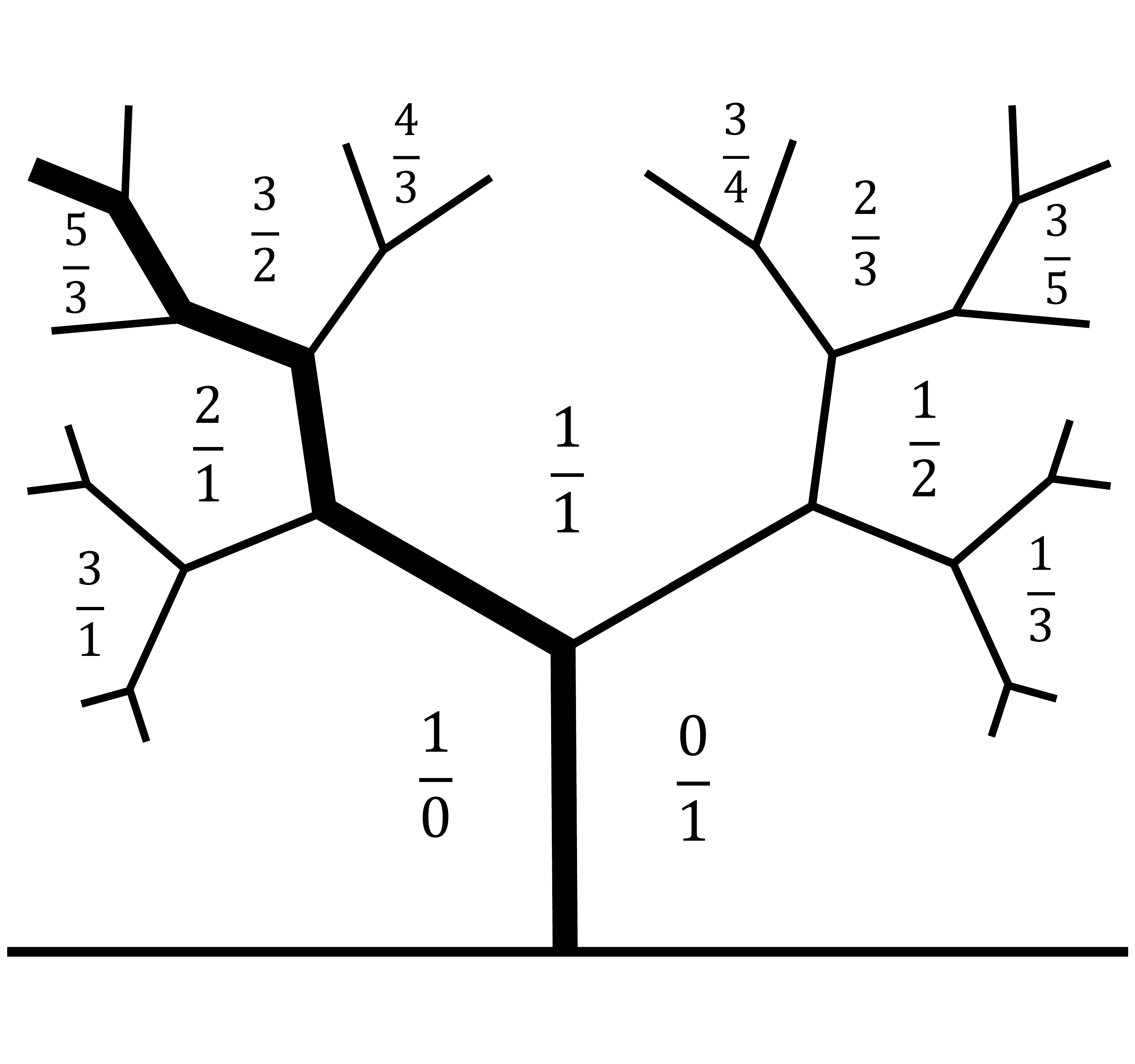}
\caption{\small Topograph of $Q=x^2+y^2$ and the corresponding positive part of the Farey tree with marked Fibonacci path.} \label{fig:golden}
\end{center}
\end{figure}

The path $\gamma_\xi$ is actually a geometric way to represent the continued fraction expansion of $\xi =\left[ a_0, a_1, a_2, a_3 \ldots \right]$:
it has $a_0$ left-turns on the tree, followed by $a_1$ right-turns, followed by $a_2$ left-turns, and so on
(see Figure \ref{fig:golden}, showing the Fibonacci path corresponding to the golden ratio $\xi=[1,1,1\dots]$).

Let us assume now that the form $Q$ is indefinite and does not represent zero, meaning that
$Q(x,y) \neq 0$ for all $(x,y)\in \mathbb Z^2\setminus (0,0).$ 

\begin{figure}[h]
\begin{center}
\includegraphics[height=95mm]{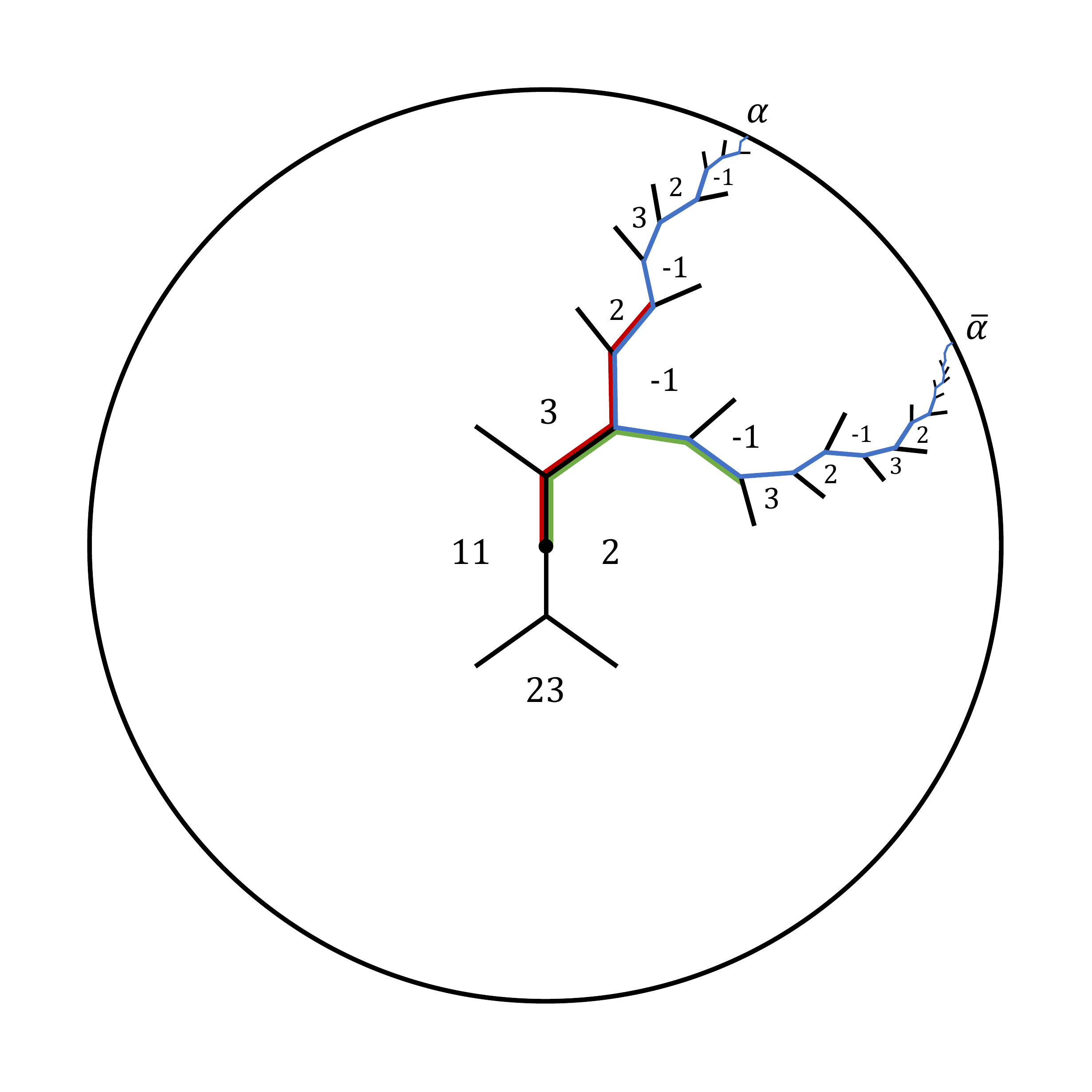} 
\caption{\small Paths to the roots $\alpha$ and $\bar\alpha$ and periodic Conway river for $Q=11x^2-10xy+2y^2.$} \label{fig:river-17}
\end{center}
\end{figure}

In this case the same arguments as in the integer case (Conway 1997) show that on the topograph of $Q$ positive and negative values are separated by an infinite path
which we call the {\it Conway river}. In the integer case we explained in (Spalding and Veselov 2017) how the Conway river is related to the continued fraction expansion of the roots $\alpha, \bar\alpha$ of the quadratic equation $Q(x,1)=0$
(see Figure \ref{fig:river-17}).

Note that in the general case of binary quadratic forms with real coefficients, we do not have periodicity of the river anymore.

\section{Proof of the theorem}

Consider an indefinite quadratic form $Q(x,y)=ax^2+hxy+by^2$ and factorise it as a product of linear forms
$$
Q(x,y)=b(y-\alpha x)(y-\beta x)
$$
with irrational $\alpha, \beta$, assuming without loss of generality that $\alpha>0$ and $\beta<0.$

Let $P=A_0$ be a corner of the Arnold sail of the corresponding pair of lines $y=\alpha x$ and $y=\beta x$. Choose a new basis in the lattice with $e_1=OA_0$ and $e_2$ being a primitive vector along the edge $A_0A_1$ of the Arnold sail. From Klein's construction it follows that this indeed a basis.

\begin{figure}[h]
\begin{center}
\includegraphics[height=90mm]{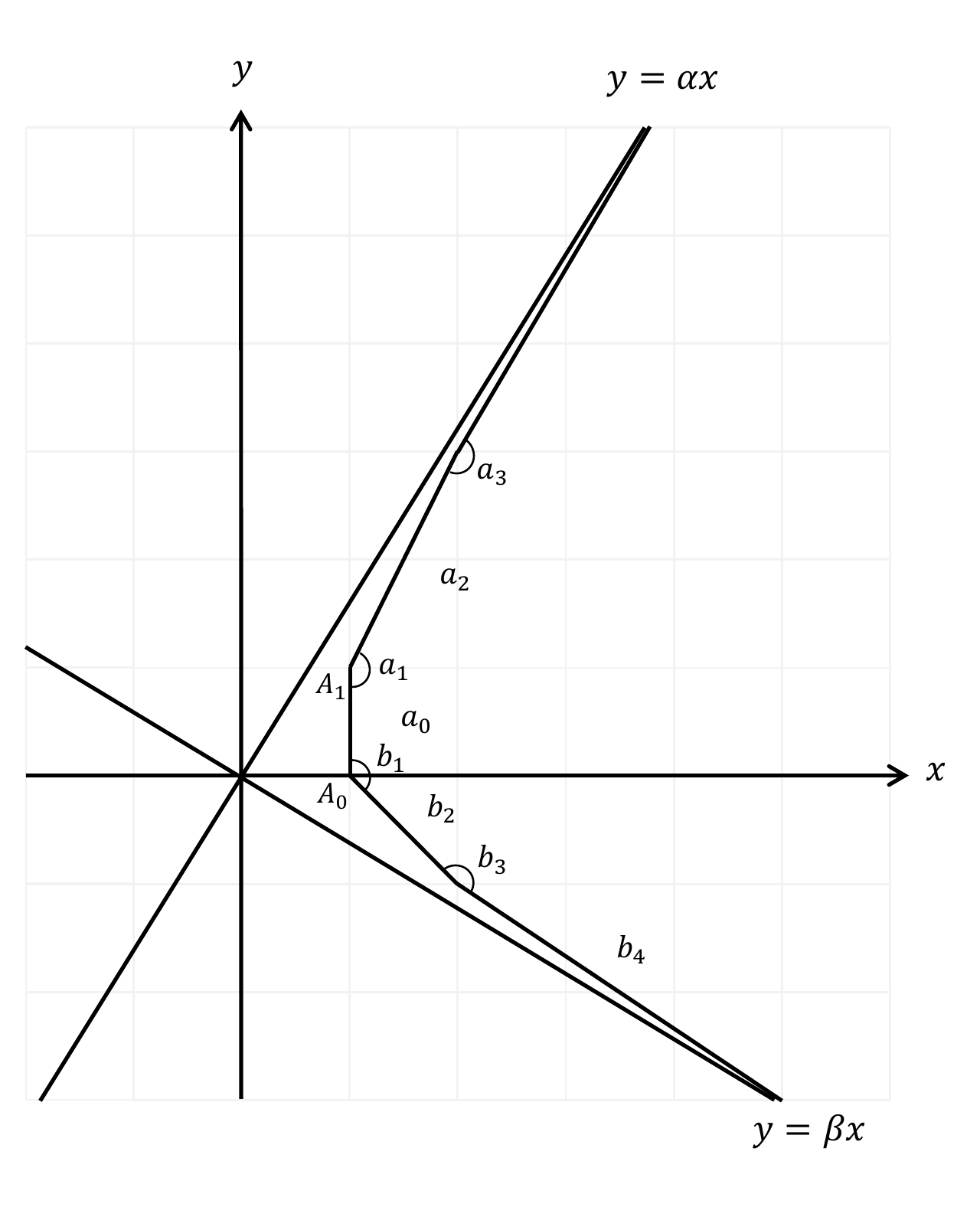} 
\caption{\small Arnold sail in a special basis.} \label{fig:normal}
\end{center}
\end{figure}

In the new coordinate system the corresponding $\alpha>1$ and $0>\beta>-1,$ and we have the situation shown in Figure \ref{fig:normal} justified by the following lemma (see also (Markov 1879) and Definition 1.1 from (Karpenkov 2018)).

\begin{Lemma}
The LLS sequence of the Arnold sail of a pair of lines $y=\alpha x$ and $y=\beta x$ with $\alpha>1$ and $0>\beta>-1$ is 
\begin{equation}
\label{llseq}
\dots, b_4, b_3, b_2, b_1, a_0, a_1, a_2, a_3, \dots,
\end{equation}
where $a_i$ and $b_j$ are given by the continued fraction expansions
\begin{equation}
\label{llseqab}
\alpha=[a_0, a_1, a_2, a_3, \dots], \quad -\beta=[0, b_1, b_2, b_3, b_4, \dots].
\end{equation}
\end{Lemma}

The proof follows directly from Klein's construction and the results of Karpenkov (2013) (see Ch. 3 and 7, in particular, Prop. 7.5).

Let us now look at the corresponding Conway river. Since $\alpha\beta<0$ this means that $Q(1,0)=b\alpha\beta$ and $Q(0,1)=b$ have different signs, so our initial position is already on the Conway river.

We know that the Conway river is the unique path on the Farey tree connecting the points $\alpha$ and $\beta$ on the boundary, and thus is the union of two paths $\gamma_\alpha$ and $\gamma_\beta$. Combining this with the description of the Farey paths in terms of continued fractions (see above), we conclude that the sequence (\ref{llseq}) determines the sequence of the river's left and right turns. 

Now let's prove that this determines the river uniquely modulo the action of $PGL(2, \mathbb Z)$, which is the symmetry group of the binary tree embedded in the plane.
Indeed, we have a well-known isomorphism $$PSL(2, \mathbb Z)=\mathbb Z_2*\mathbb Z_3.$$ This allows us to define the action of $PSL(2, \mathbb Z)$ on the tree with generators of $\mathbb Z_2$ and $\mathbb Z_3$ acting as rotations by $\pi$ about the edge centre and by $2\pi/3$ about the vertex. The element $diag\,(-1,1) \in GL(2,\mathbb Z)$ acts by the natural reflection and changes the orientation.

Using this action one can transform any directed edge to any other. After that the sequence of left- and right-turns determines the river uniquely. The left-right symmetry corresponds to the reflection. This proves our theorem.

\section{Concluding remarks}

Arnold sails can be defined in a much more general situation, in particular, for cubic binary forms and multidimensional simplicial cones. This is related to the geometric theory of multidimensional continued fractions, also going back to Klein (see (Karpenkov 2013) for the details). It is an interesting question as to whether there is an analogue of the Conway topograph here.

Another interesting question is to study the growth of values of the real binary quadratic forms along the paths on the Conway topograph, similar to the integer case considered in (Spalding and Veselov 2017). Note that in the real case the values of the form along the Conway river may approach zero (see e.g. Kleinbock and Weiss (2015)), so the situation here is more subtle.

\section{Acknowledgements} We are very grateful to Oleg Karpenkov for explaining us his important results about LLS sequences, and to Nikolai Andreev for helpful discussions.

The work of K.S. was supported by the EPSRC as part of PhD study at Loughborough.

\end{document}